\documentclass[12pt]{article}
\usepackage[english]{babel}
\usepackage{hyperref}
\usepackage{amsfonts}
\date{}
\newtheorem{theorem}{\bf Theorem}[section]
\newtheorem{remark}{\it Remark}[section]

\newcommand{\qed}{\ \hfill\hspace*{\fill}
$\vbox{\hrule\hbox{\vrule height1.3ex\hskip1.3ex\vrule}\hrule}$
\hss\vskip\topsep\relax}

\makeatletter \@addtoreset{equation}{section} \makeatother
\begin{document}

\title{EXPLICIT FORMULAS FOR THE MOMENTS OF THE SOJOURN TIME
IN THE M/G/1 PROCESSOR SHARING QUEUE WITH PERMANENT JOBS}

\author{S.F.Yashkov\thanks{This research  was supported in part
by Grant no. Sci.Sch.--934.2003.1 (Head R.A.Minlos).}\\
{\small Institute for Information Transmission Problems}\\
{\small 19,Bolshoi Karetny Lane,
\small 101447 Moscow GSP--4, Russia.}\\
{\small E-mail: yashkov@iitp.ru}\\
{\small fax: (7--095)209--05--79}}
\maketitle


\begin{abstract}
We give some representation about recent achievements in analysis
of the M/G/1 queue with egalitarian processor sharing discipline
(EPS). The new formulas are derived for the $j$-th moments ($j
\in {\mathbb N}$) of the (conditional) stationary sojourn time in
the M/G/1---EPS queue with $K$ ($K \in 0 \cup {\mathbb N}$)
permanent jobs of infinite size. We discuss also how
to simplify the computations of the moments.
\end{abstract}



\section{Introduction}

Processor sharing queues, made very attractive models by the
works of Kleinrock  \cite{Kl67}, \cite{Kl76} and Yash\-kov
\cite{Y83}, \cite{Y87}, 
play a central role in
queueing theory. These models were originally proposed to analyze
the performance of scheduling algorithms in time--sharing computer
systems, and continue to find new applications which pose
interesting mathematical problems. Over the past few years, the
processor sharing paradigm has emerged as a powerful concept for
modeling of Web servers, in particular, for evaluating the flow--level
performance of end--to--end flow control mechanisms like
Transmission Control Protocol (TCP) in Internet. 

The mathematical analysis of processor sharing queues has
resulted in many insightful results. Yet, a number of challenging
problems remains to be explored. The main goal of this paper is to
gain understanding  the problem of the moments of the stationary
sojourn time in the M/G/1 queue with egalitarian processor
sharing (EPS), and to derive the formulas for the $j$-th moments
($j \in {\mathbb N}$) of the (conditional) sojourn time in the
M/G/1---EPS queue with $K$ ($K \in 0 \cup {\mathbb N}$) permanent
jobs of infinite size. Our results complement and develop the
corresponding sections of the paper by Yashkova and Yashkov \cite{YY03}.

An idea of the EPS discipline\footnote{Under the EPS
discipline, the processor (server) is shared equally by all jobs
in the system. To put more concretely, when $1 \leq n < \infty$
jobs are present in the system, each job receives service at rate
$1/n$. In other words, all these jobs receive $1/n$ times the
rate of service which a solitary job in the processor would
receive. Jumps of the service rate occur at the instants of
arrivals and departures from the system. Therefore, the rate of
service received by a specific job fluctuates with time and,
importantly, its sojourn time depends not only on the jobs in the
processor at its time of arrival there, but also on subsequent
arrivals shorter of which can overtake a specific job. This makes
the EPS system intrinsically harder to analyze than, say, the
First Come
--- First Served (FCFS) queue.}
was introduced by Kleinrock
\cite{Kl67}  who studied only M/M/1 case as a limit of the
round--robin queue. In particular, he first showed that the mean
sojourn time conditioned on the initial job size (service
requirement) of the job is linear function of the size of the
job. For an overview of the literature on processor--sharing
queueing systems we refer to Kleinrock \cite{Kl76} (1976),
Koboyashi and Konheim \cite{KK} (1977), Jaiswal \cite{J82}
(1982), Yashkov \cite{Y87} (1987), \cite{Y90} (1990) and Yashkova
and Yashkov \cite{YY03} (2003).

The exact determination of the stationary sojourn time
distribution in the M/G/1---EPS queue was an open problem for a
long time. After puzzling researchers for 15 years, Yashkov
\cite{Y81} (1981), \cite{Y83} (1983) found an analytic solution
of this problem in terms of double Laplace tranforms (LT) (all
details contains  also his book \cite{Y89} (1989)).
Schassberger \cite{Sch} (1984) provided 
another (completely different) approach
to the exact solution by considering the EPS discipline as a
limit of the round--robin model (in discrete time). Later similar
solutions were also made 
by means of the variants of the methods from \cite {Y83} and from
\cite{Sch} (or their combinations). See, for example, van
den Berg \cite{vdB} (1990) or Whitt \cite{W98} (1998) 
(here we do have no possibility 
to discuss the contributions of other authors 
(Brandt and Brandt \cite{BB} (1998), Asare and Foster (1983),
Nunez--Queija (2000), Cheung et al. \cite{CBB}(2005), et alii)
to the closely related problems). We only mention that the
EPS queue with permanent jobs has been studied in \cite{vdB, W98, YY03, BB,
CBB} from  point of view which is different from our approach.
A telecommunication system with CPU scheduling
under SCO--UNIX can be
considered as an example of using of the EPS model with permanent jobs for
its description and predicting delays of the jobs.

In fact, our method
has turned out to be a very fruitful to derive many further
results, for example, the time-dependent queue-length  and
sojourn time distributions in this and related models (see, for
example, \cite{Y89, Y90, YY97, YY03}\footnote{Indeed, the
assumptions which required to use the steady--state (stationary)
solutions of any queueing systems are rarely satisfied in real
life. To be able really to apply queueing results in design and
analysis of technical systems, in very many cases, the obtained
results of steady state analysis are not sufficient. For example,
it is often necessary to investigate the behaviour of the queue
while it progresses towards a steady state (if and when a steady
state exists). Even the average queue length at time $t$ gives us
much more information in comparison with the stationary mean of
the number of jobs.

However,  few stochastic systems are known to have exact
time--dependent (transient) solutions for the distributions of
the processes. As a rule, such systems are the M/G/1 queues with
simpler 
disciplines (for example, FCFS,  
see, for example, Tak\'acs \cite{T}).
Besides, all 
time--dependent solutions of the queues of the type M/G/1 are
obtained in terms of double transforms (on space and time) from
which it is very hard to extract necessary information concerning
the behaviour of the system. (Moreover, much more advanced
mathematical techniques become necessary for the time--dependent
solutions in comparison with steady state analysis.) Some
exceptions  
give  variations of the M/M/1---FCFS queue for which closed--form
transient solutions are known. As a rule, the exact transient
analysis of the M/M/1---FCFS queue involves infinite sums of Bessel
functions. 
In general, explicit
exact solutions are highly unlikely for the time--dependent cases.}).
These results hold for any stability condition. Besides, the
entire transient and equilibrium behaviour of the M/G/1---EPS
queue is contained in the results mentioned, and the most (if not
all) available at present analytic solutions  (and also many new)
can be derived from them as special cases via standard arguments
(for example, by means of the Abel's/Tauber's theorems).
However, we shall not consider the transient solutions in this paper.

The rest of the paper is organized as follows. In Section 2 we
introduce some notations and describe our starting point 
represented by Theorem \ref{K}. In Section 3 we obtain some
interesting consequences of Theorem \ref{K},
some of which were 
proved earlier as self--contained theorems but now ones are
derived as special cases. 
The final section contains few closing remarks.


\section{Preliminaries}

In this section we give a short review of the  M/G/1---EPS queue
with $K$ permanent jobs (only in steady state).  For the
time-dependent results we refer to \cite{YY03} and also to
\cite{Y89, Y90}.

Jobs arrive to the single processor (server) according to a
Poisson process with the rate $\lambda > 0$. Their sizes
(reguired service times) are i.i.d. random variables with a
general distribution function $B(x)$ ($(B(0)~=~0,\; B(\infty)=1)$)
with the mean $\beta_1 < \infty$ and the Laplace--Stieltjes
transform (LST) $\beta(s)$. Let $\beta_{j}$ denote the $j$--th
moment of $B(x)$, $ \; j \in {\mathbb N}$. The service discipline is
the EPS: every job  is being served with rate $1/n$, when $n>0$
jobs are present in the system. The EPS discipline is modified by
having $K\ge 0$ extra permanent jobs with infinite sizes. The system
works in steady state. In other words, $\rho=\lambda \beta_{1} <1$
and very long time went from the instant $0$ that marks the start
of the work of our system till current time.

It is well known, due to Sakata et al. \cite{S69}, that the
stationary distribution $({\sf P}_{n})_{n \ge 0}$ of the
number of ordinary jobs in the M/G/1---EPS queue as $K=0$ is geometrically
distributed

\begin{equation}\label{a0}
\sf P_{0n}= (1-\rho)\rho^{n}, \quad n \in 0 \cup \mathbb N,
\end{equation}
where $\rho =\lambda \int_{0}^{\infty}(1-B(x))dx < 1$.
We note that $(\sf P_{0n})_{n \ge 0}$ depends on the service
time only through its mean.

For $K\ge0$ the equality (\ref{a0}) takes the form
\begin{equation}\label{a0K}
\sf P_{Kn}= (1-\rho)^{K+1}{n+K\choose K}\rho^{n}, \quad n
\in 0 \cup \mathbb N.
\end{equation}

We shall let that $V_{K}(u)$ denotes the conditional sojourn time of a job
of the size $u$ upon its arrival. This job enter into
the EPS system with $K \ge 0$ permanent jobs in steady state. Let
$v_{Kj}=\mathbb E[V_{K}(u)^{j}]$. (We shall omit the index $K$ in these
and similar notations when $K=0$.)

Define the LST of $V_{K}(u)$ by $v_{K}(r,u)=\mathbb E[{\rm e}^{-rV_{K}(u)}]$
for ${\rm Re}\;r \ge 0$ and $u\ge 0$.
Let  $\pi(r)$ 
be the LST of the busy period distribution (due to ordinary, that
is, non--permanent jobs). In other words, it is the positive root
of the well--known Tak\'acs functional equation \cite{T}

\begin{equation} \label{a1}
\pi(r) = \beta (r + \lambda - \lambda \pi(r))
\end{equation}
with the smallest absolutely value.

It is known from \cite{YY03} the following theorem
\begin{theorem}\label{K}
When $\rho < 1$, 
\begin{equation}\label{k}
v_{K}(r,u)\doteq {\mathbb E}[{\rm e}^{-rV_{K}(u)}]=v(r,u)^{K+1},
\end{equation}
where $v(r,u)$ is given by the the equality (\ref{a7}):
\begin{equation}\label{a7}
v(r,u) \doteq {\mathbb E}[{\rm
e}^{-rV(u)}]=\frac{(1-\rho){\rm e}^{-u(r+\lambda)}}
{\psi(r,u)-\tilde{a}(r,0,u)} \, .
\end{equation}
Here
\begin{equation}\label{a8}
\tilde{a}(r,0,u)=\lambda\psi(r,u)\ast\left[{\rm
e}^{-u(r+\lambda)} (1-B(u))\right]+\lambda{\rm
e}^{-u(r+\lambda)}\int_u^\infty(1-B(x))dx,
\end{equation}
where ``$\ast$" is the Stieltjes convolution sign (on variable
$u$), and $\psi(r,u)$ is the LST (with respect to $x$) of some
function $\Psi(x,u)$ of two variables (possessing the probability
density on variable $x$), which, in turn, has a Laplace transform
(LT) with respect to $u$
(argument $q$)
\begin{equation}\label{a3}
\tilde{\psi}(r,q) = \frac{q+r+ \lambda \beta (q+r+ \lambda)}
{(q+r+ \lambda)(q+ \lambda \beta (q+r+ \lambda))} \quad (r \geq
0, \, q>-\lambda \pi(r)).
\end{equation}
In (\ref{a3}), $\beta(r)=\int_{0}^{\infty}{\rm e}^{-rx}\ dB(x)$
and $\pi(r)$ (in the conditions imposed on (\ref{a3})) is
understood as the minimal solution of the functional equation
(\ref{a1}).
\end{theorem}

Thus, the function $\tilde{\psi}(r,q)$  is given in the form of
the two--dimensional transform of the function $\Psi(x,u))$
\begin{equation}
\tilde{\psi}(r,q)=\int_{0}^{\infty}\int_{0}^{\infty}{\rm
e}^{-rx-qu}d_{x}\Psi(x,u)du.
\end{equation}
In other words, $\psi(r,u)$  in equality (\ref{a7})
is the Laplace transform inversion
operator, $\psi(r,u)={\cal L}^{-1}(\tilde{\psi}(r,q))(r,u)$,
that is, the contour Bromvich integral
$$\psi(r,u)=\frac{1}{2\pi i}\int_{-i\infty +0}^{+i\infty +0}
\tilde{\psi}(r,q) {\rm e}^{qu} \ dq.
$$
\begin{remark}
Briefly, we have derived the expression for ${\mathbb E}[{\rm
e}^{-rV_{K}(u)}]$ by writing the sojourn time  as some
generalized functional on a 
branching process (like
the processes by Crump--Mode--Jagers) by means of simple
extensions of (non--trivial) arguments from \cite{Y81, Y83}. Using
the structure of the branching process, we found and solved a
system of partial differential equations (of the first order)
determining the components of a (non--trivial, too)
decomposition of $V_{K}(u)$.
It leads to ${\mathbb E}[{\rm e}^{-rV_{K}(u)}]$ (see also Remark \ref{dec}).
\end{remark}

\section{Results}

We showed in the Section 2 that the determination of the
steady--state sojourn time distribution in the queue M/G/1---EPS
with $K$ permanent jobs is simple extension of the results from
\cite{Y83}, \cite{Y81}. However, the solution contains the
Bromwich countour integrals. First we consider the case $K=0$.
Equivalent form of (\ref{a7}) (without contour integrals) is
given in the following theorem.

\begin{theorem}\label{K1}
Equivalent form of (\ref{a7}) (without the Bromwich countour
integrals) is given by
\begin{equation}\label{a9}
\frac{1}{v(r,u)}=\sum_{n=0}^{\infty}\frac{r^n}{n!}\xi_n(u),
\end{equation}
where
\begin{equation}\label{a10}
\xi_0(u)=1, \quad  \xi_n(u)= \frac{n}{(1-\rho)^n}u^{n-1} \ast
W^{(n-1)\ast}(u), \quad n=1,2, \ldots
\end{equation}
Here $W^{(n-1)\ast}(u)$ is $(n-1)$--fold convolution of the
steady--state waiting time distribution $W(u)$ in the familiar
M/G/1---FCFS system with itself $(W^{0\ast}(u)={\bf 1}(u)$,
$W^{1\ast}(u)=W(u))$, the LST of $W(u)$ is given by the
well--known Pollaczek--Khintchine formula as
\begin{equation}\label{PK}
w(q)= \frac{1-\rho}{1-\rho f(q)},
\end{equation}
where $f(q)=(1-\beta(q))/(q\beta_1)$ is the LST of the excess of $B(\cdot)$,
that is, $F(x)=\beta_1^{-1}\int_0^x(1-B(y)) dy$
($F^{0\ast}(x)={\bf 1}(x)$, the Heaviside function,
$F^{1\ast}(x)=F(x)$).
\end{theorem}

{\sf Proof.} We rewrite (\ref{a7})  in the form of Theorem 3.2
from \cite{Sch} (see also (5.5) in \cite{Y87}), namely
\begin{equation}\label{s1}
v(r,u)=\frac{(1-\rho)\delta(r,u)}{1-\rho \delta(r,u)
\left[\int_{0}^{u}\frac{dF(x)}{\delta(r,u-x)}+
(1-F(u))\right]}\quad ({\rm Re} \ r \geq 0),
\end{equation}
where
\begin{equation}\label{d3.5}
\delta(r,u)={\rm
e}^{-u(r+\lambda)}\left/\psi(r,u)\right.
\end{equation}
and $F(x)$ is
introduced in Theorem \ref{K1}. To reach our aim, it is used the
LT of $1/\delta(r,u)$ with respect to $u$ (argument $q$), which
is found from (\ref{a3}) as $\tilde{\psi}(r,q-r-\lambda), \quad
r\geq 0, \ q>r+\lambda- \lambda\pi(r)$ (cf. also the third line
on p.8 in \cite{Y87}). Now we obtain after simple algebra the
following power series expansion of the LT of the function
$1/v(r,u), \quad r\geq 0, \ u \geq 0$
\begin{eqnarray}\label{a12}
\int_0^{\infty}{\rm e}^{-qu}\frac{1}{v(r,u)}du &=&
\frac{1}{q}\left[1+\frac{1}{1-\rho} \frac{r}{q} \frac{1}
{1-\frac{1}{1-\rho} \frac{r}{q} w(q)}\right] \nonumber \\
&=& \frac{1}{q}\left[1+\sum_{n=1}^{\infty}
\left(\frac{1}{1-\rho}\frac{r}{q}\right)^n w(q)^{n-1}\right]
\end{eqnarray}
where $w(q)$ is given by (\ref{PK}). We note that
$\left|\frac{rw(q)}{(1-\rho)q}\right|<1$ as $q>r+\lambda-\lambda
\pi(r), \ \rho<~1$. Now it is easily to invert analytically (on
argument $q$) each term of the power series in $r$ (\ref{a12}).
The result is given by (\ref{a10}) whence it follows (\ref{a9}),
the right--hand side of which is the power series in $r$ with
coefficients $\xi_n(u)/n!$. \qed

The idea of such approach goes back to Heaviside. Similar results
are obtained in \cite{YY99}, \cite{ZB}. 
In fact, the form of ${\sf Var}[V(u)]$
\cite{Y81}, \cite{Y83} (see the equality (\ref{m12}) below)
stimulates a guess about the possibility of such expansion.

\begin{remark}
The formula for $W^{n \ast}(x)$ in (\ref{a10}) can be represented
in the following form
$$
W^{n \ast}(x)= (1- \rho)^{n} \sum_{k=0}^{\infty} {k + n - 1
\choose n - 1} \rho^{k}F^{k*}(x).
$$
It is done, for example, by inversion of $w(q)^{n}$, where $w(q)$
is given by (\ref{PK}).
\end{remark}

\begin{remark}\label{r3.4}
We note that the by--product of our analysis is the distribution
function $W(x)$ whose LST is given by (\ref{PK}). However, the
analysis of EPS queue gives the other quantity (corresponding to
a non--probability measure) $W^{\circ}(x)=W(x)/(1-\rho)$. The
form of the LST of $W^{\circ}(x)$ is well--known:
$w^{\circ}(q)=\sum_{n=0}^{\infty} \rho ^n f^n(q)$. Unlike $W(x)$,
$W^{\circ}(x)$ is well defined for all $\rho > 0$ and $x > 0$. It
can be shown that $W^{\circ}(x) < \infty$ for all $\rho > 0$, $x
> 0$ and for any $B(\cdot)$ (despite on the fact that,
for $\rho \geq 1$, $W^{\circ}(x) \to \infty $ as $x \to \infty$).
\end{remark}

\begin{theorem} \label{c6}
Let $v_n(u)={\mathbb E}[V(u)^n], \, n=1,2, \ldots$ Then it holds
the following recursive formula
\begin{equation}\label{a13}
v_n(u)=\sum_{i=1}^n{n \choose i}v_{n-i}(u)\xi_i(u)(-1)^{i+1}
\end{equation}
\end{theorem}

{\sf Proof.} Because $v(r,u)$ is analytical function in $r$ (in
particular, in $r=0$), we can use the Tailor series expansion of
$v(r,u)$ for small $r>0$
\begin{equation}\label{a14}
v(r,u)=1-
\frac{r}{1!}v_1(u)+\frac{r^2}{2!}v_2(u)-\frac{r^3}{3!}v_3(u)+
\ldots
\end{equation}
The product of (\ref{a14}) and (\ref{a9}) gives
\begin{eqnarray}
-\frac{r}{1!}[v_1(u)-\xi_1(u)] &+&
\frac{r^2}{2!}[v_2(u)-2v_1(u)\xi_1(u)+\xi_2(u)] \nonumber \\
&-&
\frac{r^3}{3!}[v_3(u)-3v_2(u)\xi_1(u)+3v_1(u)\xi_2(u)-\xi_3(u)]+
\ldots=0 \nonumber
\end{eqnarray}
and it leads to (\ref{a13}) after differentiating $n$ times with
respect to $r$ and setting $r=0$. \qed

In particular, the expressions for the first two moments of $V(u)$
are:
\begin{equation}\label{a15}
v_{1}(u)={\mathbb E}[V(u)]=u/(1-\rho)
\end{equation}
(this is well-known result due to Sakata et al. \cite{S69}
(1969)),
\begin{equation}\label{m12}
{\sf Var}[V(u)]=v_{2}(u)- v_{1}^{2}(u)= \frac{2}{(1-\rho)^{2}}
\int_0^u (u-x)(1-W(x)) \, dx,
\end{equation}
where  $W(x)$ is introduced in Theorem \ref{K1}, and it is
expressed as
\begin{equation}\label{m13}
W(x)= (1-\rho) \sum _{n=0}^{\infty} \rho^n F^{n*}(x)
\end{equation}
(other variables were introduced above).

The formula (\ref{a10}) implies that $\xi_1(u)={\sf E}[V(u)]$ in
(\ref{a15}). The formula  for the conditional variance
(\ref{m12}) was first obtained by Yashkov \cite{Y81}.
The standard way for the
computation of the moments is the following
\begin{equation}\label{Y0}
v_{n}(u)=\lim_{r \downarrow 0}(-1^{n})\frac{\partial^{n}v(r,u)}
{\partial r^{n}},
\quad n \in {\mathbb N}.
\end{equation}
However, the LST
$v(r,u)$ in Theorem \ref{K} is very hard to differentiate in $r$
more than once (practically almost impossible matter) since this
LST has a rather complex form due to a highly complicated form of
such constituents of (\ref{a7}) as $\tilde{a}$ and $\psi$.
Therefore ${\sf Var}[V(u)]$ is first obtained by solving
an alternative system of
differential equations (see, for example,  \cite[Chapter 2]{Y89}) which are
derived by analogy with the equations of \cite[Section 2]{Y81} or
with the equations in the proof of \cite[Theorem 4]{Y83}. These
equations are simpler forms of equations from \cite{Y81},
\cite{Y83} because ones are composed not for the LST $v(r,u)$ but
only for the second and the first moments.  Thus the formula
(\ref{a13}) for the case $n=2$ in Theorem \ref{c6} was derived
20 years earlier than the same
formula  for arbitrary integer $n$ (see \cite{YY99}, 
\cite{ZB}) 

It can be useful for asymptotic expansion of $v_n(u)$ for small
and large $u$ in the spirit of such expansion for ${\sf
Var}[V(u)]$ (see final section for some details).
Such results for ${\sf Var}[V(u)]$ were obtained
at first in \cite{Y81} (1981) (see also
\cite{Y83}).

Now we show how to extend the formulas (\ref{a15}) and
(\ref{m12}) to the case when
the M/G/1---EPS queue is modified by having $K \geq 0$ extra
permanent jobs with infinite sizes.

\begin{theorem} \label{co}
In the setting above,
\begin{equation}\label{Y1}
{\mathbb E}[V_{K}(u)]=\frac{(K+1)u}{1-\rho},
\end{equation}
\begin{equation}\label{Y2}
{\sf Var}[V_{K}(u)]=\frac{2(K+1)}{(1-\rho)^2} \int_0^u (u-x)(1-W(x)) dx,
\end{equation}
where $W(u)$  is the steady--state waiting time distribution in the
M/G/1---FCFS queue, represented by the equality (\ref{m13}).
\end{theorem}

{\sf Proof.} In our case, the equality (\ref{Y0}) takes the form
\begin{equation}\label{Y0K}
v_{K n}(u)=\lim_{r \downarrow 0}(-1^{n})\frac{\partial^{n}v_{K}(r,u)}
{\partial r^{n}},
\quad n \in {\mathbb N}.
\end{equation}
The formula (\ref{Y1}) follows directly from (\ref{k})
by means of applying (\ref{Y0K}) as $n=1$.

Taking into account (\ref{m12}), the formula (\ref{Y2})
follows also from (\ref{k})
by means of applying (\ref{Y0K}) as $n=2$ after some simple algebra. \qed

\begin{remark}\label{dec}
An alternative way to obtain (\ref{Y2}) is the following.
We can compose and solve
the system of the partial differential equations
(of the first order) which satisfy the second and the first moments
of $V_{K}(u)$. The variant of such equations is known from \cite{Y90, Y89}
as $K=0$. We point out the following fact. These equations  rely on a
decomposition of the sojourn time of the (tagged)job with the size $u$
that arrives to the EPS queue when $n$ standard jobs are present
with remaining service demands $x_1,\dots x_n$ (a key ingredient
of analysis). Denoting this
conditional sojourn time by $V_{K n}(u;x_1,\dots x_n)$, it holds
\begin{equation} \label{x5}
V_{K n}(u;x_1,\dots x_n)\stackrel{d}{=}(K+1)D(u)+\sum_{i=1}^n
\Phi(x_i,u),
\end{equation}
where all components are independent random variables.

The random variable $D(u)$ constitutes a \lq \lq main'' component
of the sojourn time: it has the distribution of the sojourn time
of a job with the size $u$ that enters into a empty
(from the standard jobs) system.
By the way, its LST is given by (\ref{d3.5}).
When the system is not empty, the $i$--th  standard job (among the jobs
which are sharing the capacity of the processor together
 with permanent jobs),
having remaining size $x_i$, \lq \lq adds'' a delay
$\Phi(x_i,u)=\Phi(x_i \wedge u,u)$ to the new job's sojourn time.
Note that $D(u)=\Phi(x_i,u)$ for $x_i \geq u$. 
Then the same chain of arguments as in \cite{Y83} can be used to
derive (\ref{k}).
\end{remark}

\section{Conclusion}

Using Theorems \ref{K} and \ref{c6}, we can easily obtain all
other moments of $V_{K}(u)$ in M/G/1---EPS queue with $K\ge 0$
permanent jobs. However, the exact expression even for the
variance of the sojourn time  (see (\ref{Y2}) involves an
integration term, making an exact computations difficult from
practical point of view. The same holds for the third, fourth,
etc. moments. This difficulty remains also in the case $K=0$. To
overcome the difficulty, it is possibly to obtain some simple
approximations for the second moments, see, for example, Villela
et al. \cite{VPR} or van den Berg \cite{vdB}. We note that there
exist also an upper (and lower) bounds for ${\sf Var}[V(u)]$.
These bounds only depends on $\rho$ and the size of job $u$. In
addition, the bounds have the attractive property of
intensitivity to $B(x)$, and the difference between the upper and
lower bounds is small, particularly, for small and moderate
values of $\rho$. These second moments tight bounds can be easily
generalized into higher moments of $V(u)$ and also to the case
$K>0$.

It are also known the asymptotic estimates of ${\sf Var}[V(u)]$ as $u \to 0$
and $u \to \infty$ \cite{Y86}. For example,
$$ Var[V(u)] \sim \frac{u^{2}\rho}{(1-\rho)^{2}} ~~{\rm as}~~u\to 0.$$
This is  some asymptotics of the sojourn time variance of a very
small jobs, and it
leads to intensitive upper bounds with special structure requiring
only knowledge of the traffic load and the job size.  Now our results
may be easily extended to the higher moments and also to the
case $K$ permanent jobs. Moreover, some preliminary analysis of
asymptotics (see \cite{Y81}) tells us about
a high accuracy of such estimates in many typical cases.
%
%

\end{document}